\documentclass[a4paper,final]{siamart171218}

\usepackage{mathtools}
\usepackage{amsmath,amssymb,amsfonts,color}
\usepackage{algorithm}
\usepackage{algorithmic}
\usepackage{color}
\usepackage{booktabs}
\usepackage[mathscr]{eucal}
\usepackage{url}
\usepackage{dsfont}
\usepackage{tabularx}
\usepackage{hyperref}
\usepackage{caption}
\usepackage{graphicx}
\usepackage{bm}
\usepackage{epstopdf}
\usepackage{xargs}
\usepackage{cleveref}
\usepackage[caption=false]{subfig}
\usepackage{algorithm,breakurl,hypcap,ifpdf,ntheorem}
\usepackage[colorinlistoftodos,prependcaption,textsize=tiny]{todonotes}

\newtheorem{remark}{Remark}
\newtheorem{assumption}{Assumption}

\newcommand{\cC}{\mathcal{C}}
\newcommand{\cF}{\mathcal{F}}
\newcommand{\R}{\mathbb{R}}

\newcommand{\cI}{\mathcal{I}}
\newcommand{\cG}{\mathcal{G}}

\newcommand{\cU}{\mathcal{U}}

\newcommand{\cN}{\mathcal{N}}
\newcommand{\cJ}{\mathcal{J}}
\newcommand{\cH}{\mathcal{H}}
\newcommand{\cL}{\mathcal{L}}
\newcommand{\cW}{\mathcal{W}}

\newcommand{\px}[1]{\partial_{{x} _{#1}}}

\newcommand{\norm}[1]{\left\lVert #1\right\rVert}
\newcommand{\bU}{\mathbf{U}}
\newcommand{\bW}{\mathbf{W}}
\newcommand{\bF}{\mathbf{F}}
\newcommand{\bH}{\mathbf{H}}
\newcommand{\bG}{\mathbf{G}}

\newcommand{\bb}{\mathbf{b}}
\newcommand{\bc}{\mathbf{c}}
\newcommand{\la}{\langle}
\newcommand{\ra}{\rangle}

\definecolor{philipp}{RGB}{205, 102, 000}

\headers{Hamilton-Jacobi-Isaacs Robust Control of Nonlinear PDEs}{D. Kalise, S. Kundu and K. Kunisch}

\title{Robust feedback control of nonlinear PDEs by numerical approximation of high-dimensional Hamilton-Jacobi-Isaacs equations}

\author{Dante~Kalise\thanks{D. Kalise is with the Department of Mathematics, Imperial College London, South Kensington Campus, SW7 2AZ London, United Kingdom (\email{dkaliseb@ic.ac.uk}).}
\and Sudeep~Kundu\thanks{S. Kundu is with the Institute of Mathematics and Scientific Computing, University of Graz, Heinrichstr. 36, A-8010 Graz, Austria(\email{sudeep.kundu@uni-graz.at}).}
\and Karl~Kunisch\thanks{K. Kunisch is with the Institute of Mathematics and Scientific Computing, University of Graz, Heinrichstr. 36, A-8010 Graz, Austria and
	Radon Institute for Computational and Applied Mathematics (RICAM), Altenbergerstra{\ss}e 69, A-4040 Linz, Austria(\email{karl.kunisch@uni-graz.at}).}}

\begin{document}
\maketitle
\begin{abstract}
We propose an approach for the synthesis of robust and optimal feedback controllers for nonlinear PDEs. Our approach considers the approximation of infinite-dimensional control systems by a pseudospectral collocation method, leading to high-dimensional nonlinear dynamics. For the reduced-order model, we construct a robust feedback control based on the $\cH_{\infty}$ control method, which requires the solution of an associated high-dimensional Hamilton-Jacobi-Isaacs nonlinear PDE. The dimensionality of the Isaacs PDE is tackled by means of a separable representation of the control system, and a polynomial approximation ansatz for the corresponding value function. Our method proves to be effective for the robust stabilization of nonlinear dynamics up to dimension $d\approx 12$. We assess the robustness and optimality features of our design over a class of nonlinear parabolic PDEs, including nonlinear advection and reaction terms. The proposed design yields a feedback controller achieving optimal stabilization and disturbance rejection properties, along with providing a modelling framework for the robust control of PDEs under parametric uncertainties.	
\end{abstract}



\section{Introduction}\label{intro}
Robustness is a feature of paramount importance in control systems design. That is, given a dynamical system onto which we act through an external control signal, we expect this signal to achieve asymptotic stabilization and to simultaneously compensate the effects of the different disturbances arising in the control loop. Perturbations may appear in the form of noisy measurements, exogenous disturbances, as well as parametric and structural uncertainties. A first step towards robustness in the control action is the synthesis of feedback laws, i.e. a control expressed as a  mapping from the current state of the system to the space of possible control actions. Feedback controls, in contrast to open-loop controls which are expressed solely as a function of time, exhibit enhanced stability properties. The design of feedback maps is a central topic in control theory, and has been thoroughly studied in its many facets, spanning ad-hoc and optimization-based approaches for linear as well and nonlinear dynamics (we refer the reader to \cite{BTB00} for a comprehensive overview).

In this paper we are particularly concerned with the design of optimal feedback controllers for nonlinear dynamics with enhanced robustness properties. The natural course of action is the use of Dynamic Programming techniques which establish a functional relation to be satisfied by the value function of the underlying control problem. In this context, once the value function associated to the control has been solved, an optimal feedback map is trivially obtained as a generalized gradient of the value function.  This approach leads to a characterization of the value function as a viscosity solution of a nonlinear Hamilton-Jacobi-Bellman (HJB) or Hamilton-Jacobi-Isaacs (HJI) type Partial Differential Equations (PDEs), depending on whether additional robustness properties are required. In the control engineering literature these approaches are known as the $\cH_2$ control synthesis whenever an optimal feedback is sought through the solution of a HJB equation, and as $\cH_{\infty}$ control design when an optimal robust feedback is computed upon the HJI equation or a related variation. The interplay between Hamilton-Jacobi-Bellman and Isaacs' PDEs and $\cH_2$ and $\cH_{\infty}$ optimal control theory have been thoroughly discussed in \cite{DG89, V17, BCD97, S92}. We specifically point the reader to \cite{S96, S99} for a study of the link between Isaacs' equation, differential game theory, and the computation of robust feedback controls via nonlinear $\cH_{\infty}$ synthesis. We recall that in the case of linear dynamical systems with no constraints on the control action and a quadratic cost, the Dynamic Programming approach reduces to the computation of Algebraic Riccati Equations, a topic which has been extensively discussed in the mathematical and engineering literature, see for instance \cite{BPDM07}. Here, our focus is on nonlinear dynamical systems, and therefore our discussion will be centered around the numerical approximation of HJB and HJI equations (HJ PDEs in short).

The numerical approximation of Hamilton-Jacobi type equations is a topic dating back to the seminal work by Crandall and Lions \cite{CL84}, including a wide range of discretization strategies such as finite differences \cite{CDI84}, finite element methods \cite{GR85}, level-set methods \cite{OF03}, domain decomposition algorithms \cite{CCFP12,F16}, filtered schemes \cite{FS} and most notably, semi-Lagrangian schemes \cite{AFK15,FF14}. In particular, the numerical approximation of the Hamilton-Jacobi-Isaacs equation has been studied in \cite{BM98, BMW99, BHT11,F06,CCF11}. For a comprehensive overview of the different approximation techniques available for HJ PDEs we refer the reader to  \cite{FF16}. The aforementioned techniques have proven to overcome the difficulties associated to the nonlinear character of the HJ PDEs. However, since they are  grid-based discretization schemes they suffer from the curse of dimensionality and become overwhelmingly expensive for problems with dimension higher than 5. This is particularly challenging in the context of nonlinear optimal control, as the dimension of the associated HJ PDE is determined by the dimension of the state space of the control problem. A partial remedy to this problem is the  coupling of model reduction techniques for the dynamics with a grid-based discretization of the HJ equation, an alternative which has been applied in \cite{AF13,AFV17, KK14, KVX04}. Another alternative is the use of approximate dynamic programming techniques in the context of reinforcement learning, see e.g. \cite{B19}. However, the design of numerical methods for the solution of high-dimensional HJ PDEs remains a daunting task. Along this direction, some encouraging results have been obtained over the last years in connection with the use of sparse grids \cite{BGGK13,GK16}, causality-free methods \cite{KW16,CDOY,YD18}, machine learning techniques \cite{SS18,RZ19}, tensor calculus \cite{YPL}, graph-tree structures \cite{AFS18}, Taylor series expansions \cite{BKP19}, and polynomial approximation  \cite{KK18}. In this latter work, we develop a numerical scheme based on a high-dimensional polynomial ansatz for the value function coupled with a Newton-type (policy iteration) method for the solution of the Galerkin residual equation associated to the HJB equation. This method allows to solve HJB PDEs up to dimension 14, becoming an effective tool for the control of dynamical systems arising from the spectral semi-discretization of nonlinear parabolic PDEs.

In this paper, our focus is on the construction of feedback laws for nonlinear parabolic PDEs which are robust with respect uncertainties and structural perturbations. These feedback laws are obtained on the basis of approximations to the  Isaacs equation. We are also particulary interested in differences between the HJB and HJI based feedback controls. As a first step towards this goal, we consider a semi-discretization in space of the infinite-dimensional control system by means of a pseudospectral collocation method \cite{qvbook,os06}. While this approach dramatically reduces the dimensionality of the resulting state space, we are still left with a nonlinear dynamical system of prohibitively large for grid-based schemes. Therefore, a fundamental step in our work is the construction of a numerical scheme for the approximation of the high-dimensional Isaacs equation.  For this, we extend our previous work \cite{KK18} to the HJI PDE by considering a high-dimensional polynomial approximation for the value function, together with the use of a bilevel policy iteration algorithm for the nonlinear Galerkin residual equation.
The final outcome of our method is an optimal feedback map, constructed by the solution of  the HJI equation, with enhanced robustness properties with respect to additive disturbances in the control loop and/or parametric uncertainties in the dynamics. The method is successfully applied to the control of nonlinear parabolic equations such as the viscous Burgers' equation and the Newell-Whithead equation under different uncertainty scenarios. To our knowledge, this is the first attempt to design optimal robust feedback controls for nonlinear PDEs by directly approximating the high-dimensional HJI equation rather than relying on the application of the classical linear $\cH_{\infty}$ method applied over the linearized dynamics.

The rest of the paper is structured as follows. In Section 2 we state the robust optimal feedback control problem for nonlinear dynamics and the associated Hamilton-Jacobi-Isaacs equation. In Section 3 we introduce an approximation scheme for the HJI PDE, and Section 4 is devoted to specific aspects of its implementation for high-dimensional dynamics. Section 5 concludes our work with a thorough computational assessment of the proposed methodology over a class of nonlinear parabolic PDEs.

\section{Preliminaries}
In the following, we discuss the synthesis of optimal feedback controllers for nonlinear dynamics based on the dynamic programming principle and more precisely, on the solution of static Hamilton-Jacobi-Bellman and Hamilton-Jacobi-Isaacs equations. In the control literature these methods are referred to as the $\cH_2/\cH_{\infty}$ control syntheses, and we briefly revisit them.
\subsection{Nonlinear $\cH_2$ control and the Hamilton-Jacobi-Bellman equation}
We consider the following infinite horizon optimal control problem:

\[\underset{u(\cdot)\in\cU}{\min}\;\cJ(u(\cdot);x):=\int\limits_0^\infty \ell(y(t))+\|u(t)\|_R^2\, dt
\]
subject to the nonlinear dynamical constraint
\begin{equation}\label{dynamics}
\dot y(t)= f(y(t))+g(y)u(t)\,,\quad y(0)=x,
\end{equation}
where we denote the state $y(t)=(y_1(t),\ldots,y_d(t))^t\in\R^d$, the control $u(\cdot)\in\cU$, with $\cU=\{u(t):\, \R_+\rightarrow U\subset\R^m\}$, the state running cost $\ell(y)\geq0$, and  the quadratic control penalization term $\|u\|_R^2=u^tRu$ with $R\in\R^{m\times m},\,R>0$. Furthermore, we assume the running cost $l(y)$, the system dynamics $f(y):\R^d\rightarrow\R^d$, and {$g(y):\R^d\rightarrow\R^{d\times m}$} to be $\cC^1(\R^d)$, $f(0)=0$ and $\ell(0)=0$. Our focus is therefore on asymptotic stabilization to the origin. By the application of the Dynamic Programming Principle, it is well-known that the optimal value function
\[V(x)=\underset{u(\cdot)\in\cU}{\inf} J(u(\cdot);x)\]
characterizing the solution of this infinite horizon control problem is the unique viscosity solution of the Hamilton-Jacobi-Bellman equation
\begin{equation}\label{hjb}
\underset{u\in U}{min}\{ DV(x)^t(f(x)+g(x) u)+ \ell(x)+\|u\|_R^2\}=0\,,\quad V(0)=0\,,
\end{equation}
with $D V(x)=(\px{1}V,\ldots,\px{d}V)^t$. We restrict our analysis to the unconstrained case $U\equiv \R^m$ where  the minimizer $u^*$ of \eqref{hjb} is given explicitly by
\begin{equation}\label{optc}
u^*(x)=\underset{u\in U}{argmin}\{ DV(x)^t(f(x)+g u)+ \ell(x)+\|u\|_R^2\}=-\frac{1}{2}R^{-1} g(x)^t DV(x)\,.
\end{equation}
By inserting \eqref{optc} into \eqref{hjb}, we obtain the HJB equation
\begin{equation}\label{hjb2}
DV(x)^t f(x)-\frac{1}{4}DV(x)^tg(x)R^{-1}g(x)^tDV(x)+\ell(x)=0\,,
\end{equation}
to be understood in the classical sense.
\begin{remark}
Under the ansatz $\ell(x)=x^tQx$, $f(x)=Ax$ and $V(x)=x^t\Pi x$ with with $Q,A,$ and $\Pi \in\R^{d\times d}$,  equation \eqref{hjb2} becomes the Algebraic Riccati Equation associated to infinite horizon linear-quadratic control. The motivation behind our work is to depart from this framework by considering a nonlinear representation of $f(x)$ and a high-order polynomial expansion for $V(x)$, in order to reduce the curse of dimensionality which arises when attempting to solve \eqref{hjb2} over high-dimensional spaces. 
\end{remark}
It is always worth to stress that by solving the HJB equation \eqref{hjb2}, the optimal feedback map $u^*(x)$ is obtained as a by-product via \eqref{optc}. For a real-time realization of the optimal control signal $u^*(t)$ in \eqref{dynamics}, the optimal feedback map is stored and evaluated online at the current state of the dynamics $y(t)$.
\subsection{Nonlinear $\cH_{\infty}$ control and the Hamilton-Jacobi-Isaacs equation}
A variation of the $\cH_2$ control synthesis is generated by considering the system dynamics
\begin{equation}\label{dynamicsp}
\dot y(t)= f(y(t))+g(y)u(t)+h(y)w(t)\,,\quad y(0)=x,
\end{equation}
where an additional disturbance signal $w(\cdot)\in\cW$, with $\cW=\{w(t):\, \R_+\rightarrow W\subset\R^p\}$ enters the system through {$h(y):\R^d\rightarrow\R^{d\times p}$}. We assume that $y=0$ is an equilibrium of the system for $u=w=0$. In this case, the $\cH_{\infty}$ control goal is to achieve both internal stability of the closed-loop dynamics and disturbance attenuation. This translates into considering the design of a feedback law $u=u(x)$ such that for given $\gamma>0$, and for all $T\geq 0$ and $w\in \cL_2(0,T)$
	\begin{equation}\label{eq:kk10}
	\int\limits_{0}^{T}\ell(x)+\|u(t)\|^2_R\, dt\leq\gamma^2\int\limits_{0}^{T}\norm{w(t)}^2_{P}\;dt,
	\end{equation}
	where $\norm{w(t)}^2_{P}=w^tPw$ with $P\in\R^{p\times p},\, P>0$.

We say that system \eqref{dynamicsp} has ${\cL}_2$-gain smaller or equal than $\gamma$,
	if \eqref{eq:kk10} holds. Subsequently we are interested in finding a
	minimum value   $\gamma^*$, such that for all $\gamma\geq \gamma^*$ it is possible to find an asymptotically stabilizing feedback law. The calculation of $\gamma^*$, the so-called $\cH_{\infty}$-norm of the system, is a challenging problem  in its own right. Here we shall obtain an estimate by means of a bisection algorithm \cite{BBK89}.

Proceeding analogously as in the $\cH_@$ control synthesis, we follow the application of the Dynamic Programming Principle. For  a given $\gamma \geq \gamma^*$, the $\cH_{\infty}$ synthesis is based on the solution of the Hamilton-Jacobi-Isaacs equation
\begin{equation}\label{hji}
\underset{u\in U}{min}\,\,\underset{w\in W}{max}\{ DV(x)^t(f(x)+g(x) u+h(x)w)+ \ell(x)+\|u\|_R^2-\gamma^2\|w\|_P^2\}=0\,.
\end{equation}
Similarly  as in the Hamilton-Jacobi equation \eqref{hjb}, given an unconstrained disturbance $\cW=\R^p$ we synthesize a stabilizing feedback law $u_{\gamma}$ and its associated disturbance $w_{\gamma}$ by setting
\begin{align}
u_\gamma(x)&=-\frac12 R^{-1}g(x)^tD V_\gamma(x)\,,\label{uopt}\\
w_\gamma(x)&=\frac{1}{2\gamma^2} P^{-1}h(x)^tD V_\gamma(x)\,,\label{wopt}
\end{align}
where $ V_\gamma(x)$ solves the Hamilton-Jacobi-Isaacs equation
\begin{align}\label{hjie}
DV_\gamma(x)^tf(x)+ \frac14 DV_\gamma(x)^tQ(x)DV_\gamma(x)+l(x)=0 \,,
\end{align}
with
\begin{align*}
Q(x)=\frac{1}{2\gamma^2}h(x)P^{-1}h(x)^t-g(x)R^{-1}g(x)^t.
\end{align*}
For further details see e.g Van Der Schaft \cite{V92}. When there is no confusion, we denote $V_\gamma(x)$ by $V$. Note that if the disturbance attenuation is neglected by taking $\gamma\to\infty$, \eqref{hjie} becomes the HJB equation related to $\cH_2$ control.

\subsection{Separability assumptions and the control of semidiscretized PDEs}
In the following sections, we will discuss the approximation of the aforedescribed HJB and HJI equations linked to $\cH_2$ and $\cH_{\infty}$ nonlinear control synthesis, respectively. In particular, we will focus on  the approximate solution of such PDEs when the dimension of the state space of the dynamics to be controlled is large ($d>10$). As it will become clear later, our numerical scheme is based on the use of globally defined polynomials for approximating the value function $V(x)$  along with a Galerkin residual equation, and we will be faced against the computation of high-dimensional integrals. In order to mitigate the computational difficulties associated with the calculation of high-dimensional integrals, the following assumption will be fundamental in the construction of our control algorithm.
\begin{assumption}\label{as:1}
The  free dynamics $f(x):\R^d\to\R^d, f(x):=(f_1(x),\ldots,f_d(x))^t$ are separable in every coordinate $f_i(x)$
\begin{equation}
f_i(x)=\sum_{j=1}^{n_{f}}\prod_{k=1}^{d}\cF_{(i,j,k)}(x_k)\,,
\end{equation}
where $\cF(x):\R^d\rightarrow\R^{d\times n_f\times d}$ is a tensor-valued  function. We shall also assume a similar separable structure for $g(x)$, $h(x)$, and $\ell(x)$.
\end{assumption}
Let us note that while this assumption includes a large class of nonlinear optimal control problems, it does not cover dynamics such as those arising in particle systems with interactions governed by distances between agents, e.g. $f(x)=f(\|x_i-x_j\|)$ (see \cite{BFK15} and references therein for control-related examples).

Here, we are interested in the application of nonlinear optimal control methods for the robust stabilization of nonlinear parabolic PDEs of the form
\begin{equation}
\partial_{t}X(\xi,t) =\cL X(\xi,t)+\cN(X(\xi,t),\nabla X(\xi,t))+\mathcal{G}(\xi,X)u(t)+\mathcal{H}(\xi,X)w(t)\,,
\end{equation}
where $\cL$ is a linear elliptic operator, $\cN$ is a nonlinearity, while $\cG$ and $\cH$ are control and disturbance operators, respectively. At this point two natural questions arise: how can we apply finite-dimensional control techniques to infinite-dimensional dynamics, and whether it is plausible to assume the separability of the dynamics. The answer to both questions rely on the semi-discretization of the dynamics in space and the use of the so-called method of lines \cite{SG09}. In fact, under the Galerkin ansatz
\begin{equation}
X(\xi,t)\approx \sum\limits_{i=1}^d X_i(t)\Psi_i(\xi)
\end{equation}
with $\{\Psi_i\}_{i=1}^{\infty}$ a suitable set of basis functions in space such as local spline approximants, spectral elements or global orthogonal polynomials, the infinite-dimensional dynamics are approximated by a $d$-dimensional dynamical system of the form
\begin{equation}
\dot X(t) =A X(t)+N(X(t))+G(X)u(t)+H(X)w(t)\,,
\end{equation}
for $X(t)=(X_1(t),\ldots,X_d(t))$ and where every term is corresponds to a separable expression. The linear operator $A$ is naturally written in separable form, and for instance, in the case of a polynomial nonlinearity $\cN(X)=p(x)$ the use of collocation points in space leads to a finite-dimensional realization of the nonlinearity consisting on fully decoupled pointwise evaluations of the polynomial at each collocation point. While the use of finite elements/splines also leads to separable representations, in this work we opt for the use of spectral/pseudospectral collocation methods as the dimension $d$ of the finite-dimensional dynamics must be kept moderate. Once a suitable separable finite-dimensional state space representation of the system has been obtained, both $\cH_2$ and $\cH_{\infty}$ control synthesis methods can be applied, amounting to the solution of HJB/HJI equations over a $d$-dimensional domain. In the following sections we shall focus on the construction of high-dimensional HJI equations arising in $\cH_{\infty}$  control. The details concerning the construction of an analogous method for HJB equation in $\cH_2$ control can be found in \cite{KK18}.

\section{Iterative solution of the Hamilton-Jacobi-Isaacs equation} This section addresses the construction of a numerical scheme for the solution of the HJI equation
\begin{align}\label{hji3}
DV_\gamma(x)^tf(x)+ \frac14 DV_\gamma(x)^t\left(\frac{1}{2\gamma^2}h(x)P^{-1}h(x)^t-g(x)R^{-1}g(x)^t\right)DV_\gamma(x)+\ell(x)=0 \,,
\end{align}
which corresponds to the unconstrained version $\cU=\R^m, \cW=\R^p$ of the Isaacs' equation
\begin{equation}\label{hjimm}
\underset{u\in U}{min}\,\,\underset{w\in W}{max}\{ DV(x)^t(f(x)+g(x) u+h(x)w)+ \ell(x)+\|u\|_R^2-\gamma^2\|w\|_P^2\}=0\,.
\end{equation}
Both versions of the HJI equation share the same computational difficulties: they correspond to nonlinear PDEs  which need to be solved over a high-dimensional space.
At first glance, a natural idea to address the quadratic nonlinearity on $DV$ arising in \eqref{hji3} is to utilize a Newton iteration over the continuous PDE or its discretization. In fact, in the context of the Hamilton-Jacobi-Bellman equation
\begin{align}
DV_\gamma(x)^tf(x)- \frac14 DV_\gamma(x)^tg(x)R^{-1}g(x)^t DV_\gamma(x)+\ell(x)=0 \,,
\end{align}
the application of Newton's method is equivalent to the well-known \textsl{policy iteration} or Howards' algorithm \cite{BMZ09} for
\begin{equation}\label{eq:hjb3}
\underset{u\in U}{min}\,\{ DV(x)^t(f(x)+g(x) u)+ \ell(x)+\|u\|_R^2\}=0\,,
\end{equation}
which we recall in \cref{alg:polit}.
\begin{algorithm}[h!]
	\begin{algorithmic}
\STATE{{\bf Input}: $u^{(0)}$ be an asymptotically stabilizing control law for the unperturbed dynamics \eqref{dynamics}, a tolerance $\epsilon$.}
\STATE{\bf While $\|u^{(i+1)}-u^{(i)}\|\geq\epsilon$, \;}
		\STATE{\bf Solve for $V^{(i)}(x)$ (policy evaluation step):
			\begin{align*}
			DV^{(i)}(x)^t\big(f+g(x) u^{(i)}(x)\big)+ \ell(x)+\|u^{(i)}(x)\|^2_R=0\,,
			\end{align*}}
\vspace{-3mm}
			\STATE{\bf Update the control (policy update step):
\vspace{-3mm}
				\[
				u^{(i+1)}(x)=-\frac12 R^{-1}g(x)^tDV^{(i)}(x)\,.
				\]	
End}
		\caption{Howards' Algorithm/Continuous Policy Iteration}\label{alg:polit}
	\end{algorithmic}
\end{algorithm}
This algorithm takes as an input an asymptotically stabilizing controller for the nonlinear dynamics \eqref{dynamics}, and iterates in two steps: a \textsl{policy evaluation step} in which the current feedback $u^{(i)}(x)$ is inserted in the HJB equation \eqref{eq:hjb3} to solve for $V^{(i)}(x)$ avoiding the nonlinearity associated to the minimization, and a \textsl{policy update step} in which the formula $u^{(i+1)}(x)=-\frac12 R^{-1}g(x)^tDV^{(i)}(x)$ updates the control. Under the assumption that the initial controller $u^0(x)$ is an asymptotically stabilizing feedback, it has been proven in \cite{BM98} that \cref{alg:polit} generates a sequence converging to the optimal feedback control $u^{(\infty)}(x)$, and its associated value function $V^{(\infty)}(x)$ which solves \eqref{eq:hjb3}.
This iterative procedure can be extended in a straightforward manner to HJI equations of the form \eqref{hjimm} by nesting inside every policy evaluation step, an inner iterative loop to find the optimal adversarial perturbation $w^{(i,\infty)}$ associated to $u^{(i)}$. A version of the policy iteration procedure for HJI equations is presented in \cref{alg:sga1} below.
\begin{algorithm}[h!]
	\begin{algorithmic}
				\STATE{{\bf Input}: $u_{\gamma}^{(0)}(x)$ be an asymptotically stabilizing control law for the dynamics \eqref{dynamics} with $w=0$, a tolerance $\epsilon$, and $\gamma\geq\gamma^*$.}
				\STATE{\bf While $\|u_{\gamma}^{(i+1)}-u_{\gamma}^{(i)}\|\geq\epsilon$,
				\STATE{\bf Set $w^{(i,0)}_{\gamma}(x)\equiv 0$, $j=0$}
\STATE{\bf While $\|w^{(i,j+1)}_{\gamma}-w^{(i,j)}_{\gamma}\|\geq\epsilon$,\\
			\STATE{\bf 	\hspace{1.5cm} Solve for $V^{(i,j)}_{\gamma}(x)$ (policy evaluation):
			\begin{align}\label{eq:aux1}
			    DV^{(i,j)}_{\gamma}(x)^t\big(f+gu^{(i)}_{\gamma}+hw^{(i,j)}_\gamma\big)+ \ell(x)+\|u^{(i)}_{\gamma}\|_R^2-\gamma^2\|w^{(i,j)}_{\gamma}\|^2_P=0\,,
			   \end{align}}
\vspace{-3mm}		
		\STATE{\bf \hspace{1.5cm} Update the disturbance:
			\begin{equation*}\label{eq:aux2}
			w^{(i,j+1)}_{\gamma}(x)=\frac{1}{2\gamma^2}P^{-1}h(x)^tDV^{(i,j)}_{\gamma}(x)\,,
			\end{equation*}
			}}}
\bf {End

			\STATE{Update the control:
\vspace{-3mm}
				\begin{equation*}\label{eq:aux3}
				u^{(i+1)}_{\gamma}(x)=-\frac12 R^{-1}g(x)^tDV_\gamma^{(i,\infty)}(x)\,.
				\end{equation*}
				\;}

			End}
		\caption{Continuous Policy Iteration for Hamilton-Jacobi-Isaacs equations}\label{alg:sga1}
	\end{algorithmic}
\end{algorithm}

The convergence of policy iteration-type algorithms for differential games has been discussed in  \cite{BMW99} and \cite{BM98} , among many others. Under the same assumption on the asymptotic stability of the initial guess $u_{\gamma}^{(0)}(x)$ for $w=0$, the algorithm has been shown to converge towards an optimal stabilizing feedback $u_{\gamma}^{(\infty)} (x)$ with an associated perturbation $w_{\gamma}^{(\infty, \infty)}$, and a value function $V_{\gamma}^{(\infty,\infty)}$, which is  the unique stabilizing solution of the HJI equation \eqref{hji3}.

\begin{remark}
	Both \cref{alg:polit} and \cref{alg:sga1} assume the existence of a globally stabilizing controller for initialization. It is often the case that only a locally stabilizing controller is available, i.e. there exists a bounded domain $\Omega\subset\R^d$ where the initial conditions are stabilized with this feedback law. In that case the convergence of the algorithms is only valid inside $\Omega$, and this must be taken into account when defining a domain of interest for deriving a numerical approximation.
\end{remark}

\begin{remark}
	\cref{alg:sga1} is valid for $\gamma\geq\gamma^*$. In order to estimate $\gamma^*$ for a given nonlinear problem, we proceed similarly as in the  linear $\cH_{\infty}$ problem \cite{BBK89}. We perform a bisection method over $\gamma$ until finding the smallest value for $\gamma$ for which \cref{alg:sga1} converges to a stabilizing feedback law.
\end{remark}

 Given an asymptotically stabilizing  initialization of \cref{alg:sga1}, and having fixed $u(x)$ and $w(x)$ at a policy evaluation step \eqref{eq:aux1}, we must solve the linear PDE for $V(x)$
\begin{align*}
DV(x)^t\big(f(x)+g(x)u(x)+h(x)w(x)\big)+ \ell(x)+\|u(x)\|_R^2-\gamma^2\|w(x)\|^2_P=0\,.
\end{align*}

In the following section, we discuss the construction of a numerical scheme for the solution of this equation over high-dimensional domains by means of a Galerkin formulation with global polynomial basis functions.

\section{ Solving the linear equations}

In the previous section we described how the solution to \eqref{hji} can be approximated iteratively by solving at each iteration level a high-dimensional linear equation of the form
\begin{equation}\label{eq4.1}
\begin{aligned}
\qquad&\cG(V_\gamma,DV_\gamma;u,w)=0, \quad V_\gamma(0)=0, \text{ where } \\
\qquad & \cG(q,p;u,w)=p^t(f+gu+hw) + \ell+\|u\|_R^2 -\gamma^2\|w\|_P^2
\end{aligned}
\end{equation}
for the unknown $V_\gamma=V_\gamma(x)$, given  $u,w,f,g$ and $\ell$ all functions of $x$.

\subsection{Galerkin Approximation}
To derive a Galerkin approximation to \eqref{eq4.1} let  $\{\phi_j\}_{j=1}^{\infty}\in C^\infty(\Omega,\R)$ be a complete set of basis functions in $L^2(\Omega,\R)$ so that
$$V(x)=\sum_{j=1}^{\infty}c_j\phi_j(x), \quad\text{with}\quad \phi_j(0)=0\quad \forall j,$$
and introduce $V_n(x)$  as the approximation of $V_\gamma$,  of the form
$$V_n(x)=\sum_{j=1}^{n}c_j\phi_j(x)\equiv \Phi_n\bc,$$
where $\Phi_n:=(\phi_1(x),\ldots,\phi_n(x))$ and $\bc=(c_1,\ldots,c_n)^t$.
The coefficients $c_j$ are obtained by imposing the Galerkin residual system
\begin{equation}\label{gHJIg}
\langle \cG(V_n,DV_n,u;w),\phi_i\rangle_{L^2(\Omega)}=0\,,\quad \forall \phi_i\in\Phi_n\,.
\end{equation}
We now compute the different terms involved in the approximation of  \eqref{gHJIg} once we apply the Galerkin ansatz to $V_\gamma^{(i,j)}$ appearing in \eqref{eq:aux1}.
 The value function $V_\gamma^{(i,j)}$ is approximated with the expansion
\begin{align}\label{eq3.2}
V^{(i,j)}_n(x)=\sum_{j=1}^{n}c^{(i,j)}_j\phi_j(x).
\end{align}
For each iteration the values of $w^{(i,j)}$ and $u^{(i)}$ are determined according to
\begin{align}
w^{(i,j)}(x)=\frac{1}{2\gamma^2}P^{-1}h^tDV^{(i,j-1)}_n(x),\quad
u^{(i)}(x)=-\frac 12 R^{-1}g^tDV^{(i-1,\infty)}_n(x)\label{eq3.1},
\end{align}
where $V^{(i-1,\infty)}_n$ denotes the value function obtained in the last policy evaluation step of the previous $i-$iteration.  We next specify the additive  terms of the Galerkin residual equation corresponding to \eqref{eq:aux1} with $V_\gamma^{(i,j)}$ replaced by  $V_n^{(i,j)}$, or alternatively  \eqref{gHJIg} with $u=u^{(i)}$ and $w=w^{(i,j)}.$ These expressions are in part already available in \cite{KK18}, but for the sake of coherence we provide all of them here. Below we shall write $\bc^{(i,j)}$ for $(c^{(i,j)}_1,\dots,c^{(i,j)}_n)^t$ ($i,j=0,\hdots, \infty$).  Also note, that superscripts $i,j$ refer to the iterations in the loops of \cref{alg:sga1}, whereas subscripts $i,j$ represent running indices of the Galerkin discretization:

\bigskip

\begin{enumerate}

\item[\bf 1)] \noindent$\pmb{\langle  DV_n^tf,\phi_i\rangle_{L^2(\Omega)}}$: The expansion $V_n$ leads to $DV_n^tf=\sum_{j=1}^n c_jD \phi_j^tf,$
and hence
\[\langle D V_n^tf,\phi_i\rangle_{L^2(\Omega)}=\bF_{(i,\bullet)} \bc\,\,,\quad \bF\in\R^{n\times n}\,,\quad\bF_{(i,j)}:=\langle D\phi_j^tf,\phi_i\rangle_{L^2(\Omega)}\,.\]
\item[\bf 2)] \noindent$\pmb{\la D V_n^t hw,\phi_i\ra_{L^2(\Omega)}}$: We use relation \eqref{eq3.1} to obtain
\begin{equation*}
\begin{aligned}
D V_n^t hw&=D V_n^t\left(\frac{1}{2\gamma^{2} }P^{-1} hh^tD V_n^{(i,j-1)}\right)\\
&=\frac{1}{2\gamma^{2} }P^{-1} \sum_{k=1}^nc_kD\phi_k^t\left( hh^t\sum_{r=1}^nc_r^{(i,j-1)}D\phi_r\right)^t,
\end{aligned}
\end{equation*}
such that
\begin{align*}
\la D V_n^t hw,\phi_i\ra_{L^2(\Omega)}&=\bH_{(i,\bullet)} \bc\,,\quad \bH\in\R^{n\times n}\,,\\
\bH_{(i,j)}&:=\frac{1}{2\gamma^{2} }P^{-1}\sum_{r=1}^nc_r^{(i,j-1)} \la h^tD\phi_kD\phi_r^t h,\phi_i\ra_{L^2(\Omega)}\,.
\end{align*}
\item[\bf 3)] \noindent$\pmb{\la D V_n^t gu,\phi_i\ra_{L^2(\Omega)}}$: By inserting the expansion \eqref{eq3.1}, we arrive at
\begin{equation*}
\begin{aligned}
D V_n^t gu&=D V_n^t\left(-\frac 12 R^{-1} gg^tD V_n^{(i-1,\infty)}\right)\\
&=-\frac 12 R^{-1} \sum_{k=1}^nc_kD\phi_k^t\left( gg^t\sum_{r=1}^nc_r^{(i-1,\infty)}D\phi_r\right)^t,
\end{aligned}
\end{equation*}
such that
\begin{align*}
\la D V_n^t gu,\phi_i\ra_{L^2(\Omega)}&=\bG_{(i,\bullet)} \bc\,,\quad \bG\in\R^{n\times n}\,,\\
\bG_{(i,j)}&:=-\frac 12 R^{-1}\sum_{r=1}^nc_r^{(i-1,\infty)} \la g^tD\phi_rD\phi_k^t g,\phi_i\ra_{L^2(\Omega)}\,.
\end{align*}
\item[\bf 4)] \noindent$\pmb{\langle \ell,\phi_i\rangle_{L^2(\Omega)}}$: Assuming $\ell(x)=x^t Q x$, leads to
\[\langle \ell(x),\phi_i\rangle_{L^2(\Omega)}=\langle x^t Q x,\phi_i\rangle_{L^2(\Omega)}\,,\qquad Q> 0\in \R^{d\times d}\,.\]

\item[\bf 5)] \noindent$\pmb{\langle\gamma^2\norm{w}^2_P,\phi_i\rangle_{L^2(\Omega)}}$ : Note that
$
\gamma^2\norm{w}^2_P= \frac{P^{-1}}{4\gamma^{2}}D V_n^{(i,j-1)t}hh^tD V_n^{(i,j-1)},$
hence
$\langle\gamma^2\norm{w}^2_P,\phi_i\rangle_{L^2(\Omega)}=(\bc^{(i,j-1)})^t\bW_{(i,\bullet)}\bc^{(i,j-1)},$
where
$\bW\in\R^{n\times n\times n}$ is given by
\[\bW_{(i,j,k)}:=\langle(h^tD\phi_{j})(D\phi_{k}h),\phi_{i}\rangle_{L^2(\Omega)}\,.\]

\item[\bf 6)] \noindent$\pmb{\langle\norm{u}^2_R,\phi_i\rangle_{L^2(\Omega)}}$ : Since
$
\norm{u}^2_R\equiv u^tRu=\frac 14 R^{-1}D V_n^{(i-1,\infty)t}gg^tD V_n^{(i-1,\infty)},
$
this leads to
$\langle\norm{u}^2_{R},\phi_i\rangle_{L^2(\Omega)}=(\bc^{(i-1,\infty)})^t\bU_{(i,\bullet)}\bc^{(i-1,\infty)},$
where
$\bU\in\R^{n\times n\times n}$ is given by
\[\bU_{(i,j,k)}:=\langle(g^tD\phi_{j})(D\phi_{k}g),\phi_{i}\rangle_{L^2(\Omega)}\,.\]

\end{enumerate}
After discretization, the discretized generalized Hamilton Jacobi Isaacs  equation \eqref{gHJIg} reduces to the  linear system for $\bc$
\[\left(\bF+\bH(\bc^{(i,j-1)})+\bG(\bc^{(i-1,\infty)})\right)\bc=\bb(\bU,\bW,\bc^{(i-1,\infty)},
\bc^{(i,j-1)})\,,\]
where $\bb$ is given by the expansion of $-\Big(\ell+\norm{u}^2_{R}-\gamma^2 \norm{w}^2_{P}\Big)$ ( terms {\bf 4} , {\bf 5} and {\bf 6} in the list above).

\subsection{Separable approximation} In order to tackle the computational difficulties associated to the curse of dimensionality and the calculation of high-dimensional integrals, we resort to \cref{as:1} concerning the separability of the different terms appearing in the control system. The  multi-dimensional basis  functions $\Phi_n:=(\phi_1(x),\ldots,\phi_n(x))$ used in the expansion of $V_n$ are generated by products of one-dimensional elements chosen from a  polynomial basis $\varphi_M:\R\to \R^M$. Here we choose the monomial basis $\varphi_M=(1,x,\ldots,x^M)^T$ with $M\in \mathbb{N}$  the degree of the polynomial, however it is also possible to use orthogonal polynomials. The multidimensional basis is then generated as a subset of the $d$-dimensional tensor product of one-dimensional bases such that the maximum total degree is fixed,
\[\Phi_n\equiv\left\{\phi\in\bigotimes\limits_{i=1}^d\varphi_M(x_i)\,,\;\text{and } deg(\phi)\leq M\right\}\,.\]
A canonical basis element $\phi_i$ in $\Phi_n$ is then
given by $\phi_i= \Pi_{j=1}^d x_j^{\nu_j}$ , with $\sum_{j=1}^d \le M$. The cardinality of the set of monomials in $d$ dimensions with a total degree not greater than $M$ is given by
\begin{equation}\label{totalcount} n=\sum\limits_{m=1}^M\left(\begin{array}{c}d+m-1\\m\end{array}\right)\,,
\end{equation}
which partially mitigates the curse of dimensionality affecting mesh-based discretizations where the dependence of the total number of degrees of freedom is exponential with respect to $d$. \Cref{tab:count}  shows the cardinality of $\Phi_n$ for different values of interest for $M$ and $d$.
\begin{table}[h]
	\centering
	\begin{tabular}{|c|cccc|}
		\hline
		& \multicolumn{4}{c|}{Full monomial basis}\\
		\hline
		$d$\textbackslash$M$& 2 & 4 & 6& 8\\
		\hline
		\rule{0pt}{4ex}
		6&27 & 209 &923&3002\\
		8& 44& 494&3002&12869\\
		10& 65& 1000&8007&43757\\
		12& 90& 1819&18563&125969\\
		14&  119 & 3059 & 38759& 319769\\
		\hline
	\end{tabular}
\vskip 2mm
\caption{Cardinality of the set of $d-$ dimensional monomials with total degree smaller or equal to $M$. The dependence on $d$ is combinatorial, mitigating the curse of dimensionality for low-degree polynomials.}\label{tab:count}
\end{table}

Having generated a set of separable multidmensional basis function, we proceed as in the previous section, obtaining the summands in \eqref{gHJIg} .  The integration is carried out over the hyperrectangle $\Omega=\Omega_1\times\ldots\times\Omega_d$. In \cite{KK18}, some of the terms are given in details. For completeness, we repeat them here.
\begin{enumerate}

\item[\bf 1)] \noindent$\pmb{\langle  DV_n^tf,\phi_i\rangle_{L^2(\Omega)}}$: In this case, using the expansion of $V_n$, we require the computation of
\[\langle D\phi_j^tf,\phi_i\rangle_{L^2(\Omega)}=\sum_{p=1}^d\langle f_p\px{p}\phi_j,\phi_i\rangle_{L^2(\Omega)}\,,\]
which is expanded by using the separable structure of the free dynamics
\[\langle f_p\px{p}\phi_j,\phi_i\rangle _{L^2(\Omega)}=\sum\limits_{l=1}^{n_f}\langle \left(\prod_{m=1}^d\cF(p,l,m)\right)\px{p}\phi_j,\phi_i\rangle_{L^2(\Omega)}\,,\]
where
\begin{align*}
&\langle\left(\prod_{m=1}^d\cF(p,l,m)\right)\px{p}\phi_j,\phi_i\rangle_{L^2(\Omega)}\\
=&\left(\prod_{\begin{subarray}{l}m=1\\m\neq p\end{subarray}}^d\int\limits_{\Omega_m}\cF(p,l,m)\phi_i^m\phi_j^m(x_m)\,dx_m\right)\left(\int\limits_{\Omega_p}
\cF(p,l,p)\phi_i^p\px{p}\phi_j^p(x_p)\,dx_p\right).
\end{align*}
\item[\bf 2)] \noindent$\pmb{\la D V_n^t hw,\phi_i\ra_{L^2(\Omega)}}$: This term
demands to calculate
$\langle h^t D\phi_kD\phi_j^t h,\phi_i\rangle_{L^2(\Omega)}.$
Using the  separable expansion, we can write $h_{i}(x)=\sum_{j=1}^{n_{h}}\prod_{k=1}^{d}H_{(i,j,k)}(x_k),$ where $H(x):\R^d\rightarrow \R^{d\times n_{h}\times d}$ is a tensor valued function.
Hence \[\langle h^t D\phi_kD\phi_j^t h,\phi_i\rangle_{L^2(\Omega)}=\sum_{l,m=1}^{d}\langle h_{l}\px{l}\phi_k\px{m}\phi_jh_{m},\phi_i\rangle_{L^2(\Omega)}=:h_{sep},\] with
$$h_{l}=\sum_{r=1}^{n_{h}}\prod_{p=1}^{d}H_{(l,r,p)}(x_p),\qquad h_{m}=\sum_{q=1}^{n_{h}}\prod_{p=1}^{d}H_{(m,q,p)}(x_p).$$
Therefore, we obtain
\begin{align*}
&h_{sep}=\sum_{l,m=1}^{d}\sum_{r,q=1}^{d}\Bigg\langle \prod_{p=1}^{d}H(l,r,p)\px{l}\phi_k\px{m}\phi_j\prod_{p=1}^{d}H(m,q,p),\phi_i\Bigg\rangle\notag\\
&=\sum_{l,m=1}^{d}\sum_{r,q=1}^{d}\Bigg(\prod_{\begin{subarray}{l}p=1\\p\neq l\\p\neq m\end{subarray}}^{d}\int_{\Omega_p}\!\!H(l,r,p)(x_p)\phi^p_k(x_p)\phi^p_j(x_p)H(m,q,p)(x_p)\phi^p_i(x_p) dx_p\Bigg)\notag\\
&\quad\left(\int_{\Omega_l}H(l,r,l)(x_l)\px{l}\phi^l_k(x_l)\phi^l_j(x_l)H(m,q,l)(x_l)\phi^l_i(x_l)dx_l\right)\notag\\
&\quad\left(\int_{\Omega_m}H(l,r,m)(x_m)\px{m}\phi^m_j(x_m)\phi^m_k(x_m)H(m,q,m)(x_m)\phi^m_i(x_m)dx_m\right).
\end{align*}
\item[\bf 3)] \noindent$\pmb{\la D V_n^t gu,\phi_i\ra_{L^2(\Omega)}}$: Similar to 2).
\item[\bf 4)] \noindent$\pmb{\langle \ell(x),\phi_i\rangle_{L^2(\Omega)}}$:
Here, the expression involves \[\langle \ell,\phi_i\rangle_{L^2(\Omega)}=\sum\limits_{j,k=1}^dQ_{(j,k)}\langle x_jx_k,\phi_i\rangle_{L^2(\Omega)}\,,\]
where $\langle x_jx_k,\phi_i\rangle_{L^2(\Omega)}$ can be expanded as

\[\langle x_jx_k,\phi_i\rangle_{L^2(\Omega)}=\left(\prod_{\begin{subarray}{l}p=1\\p\neq j\\p\neq k \end{subarray}}^d\int\limits_{\Omega_p}\!\!\phi_i^p(x_p)\,dx_p\right)\!\!\left(\int\limits_{\Omega_j}\!\!\phi_i^j(x_j)x_j\,dx_j\right)\!\!
\left(\int\limits_{\Omega_k}\!\!\phi_i^k(x_k)x_k\,dx_k\right)\,.\]

\item[\bf 5)] \noindent$\pmb{\langle\norm{w}^2_P,\phi_i\rangle_{L^2(\Omega)}}$ : This term requires the computation of the inner product $$\langle(h^ tD\phi_{j})P(D\phi_{k}h),\phi_{i}\rangle_{L^2(\Omega)},$$ which is expanded with the help of $\pmb{\la D V_n^t hw,\phi_i\ra_{L^2(\Omega)}}$ as
$$\langle(h^ tD\phi_{j})P(D\phi_{k}h),\phi_{i}\rangle_{L^2(\Omega)}=Ph_{sep}.$$
\item[\bf 6)] \noindent$\pmb{\langle\norm{u}^2_R,\phi_i\rangle_{L^2(\Omega)}}$: This term is similar to 5).
\end{enumerate}
\section{Numerical experiments}
In this section we present numerical experiments illustrating  that a combination of spectral techniques for the semi-discretization of a PDE control system, together with a polynomial approximation for the value function solving the associated HJI equation can mitigate the curse of dimensionality. We apply this methodology for the robust feedback design for nonlinear parabolic PDEs.

\subsection{Convergence of the polynomial approximation}
\begin{figure}[!ht]
	\centering
	\includegraphics[width=.48\textwidth]{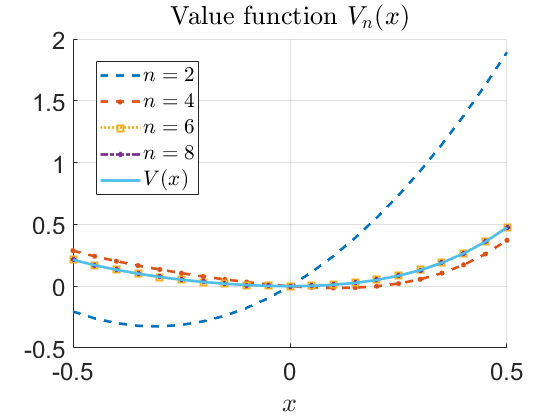}
	\includegraphics[width=.48\textwidth]{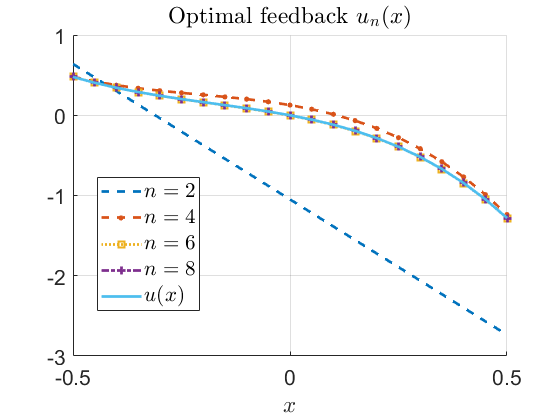}
	\caption{1D polynomial approximation of the HJI equation with nonquadratic running cost. Approximation with monomial basis. The number $n$ denotes the total number of basis functions.}\label{figconv1d}
\end{figure}

We begin by assessing the convergence of the polynomial approximation for the HJI equation in a 1D test, setting
\[f(x)=0\,,\quad g=1\,,\quad h=0.1\,, \quad l(x)=\frac{q}{4}{\left(x^2\, \mathrm{e}^{x} + 2\, x\, \mathrm{e}^{x} + 4\, x^3\right)}^2\,,\quad q=\frac{g^2}{R}-\frac{h^2}{2\gamma^2P}\,,\]
with $R=P=1$, and $\gamma=2$. In this case the exact solution of the HJI equation  is given by
\[ V(x)=x^4+x^2e^x\,.\]
We initialize \cref{alg:sga1} with $u^0=-0.01x$ and a threshold value $\epsilon=1\times 10^{-6}$. The relative error for an $n-$degree approximation of $V(x)$, denoted by $V_n(x)$, is defined as
\[\text{error } V_n:=\frac{\|V_n(x)-V(x)\|_{L^2(\Omega)}}{\|V(x)\|_{L^2(\Omega)}},\]
with $\Omega=[-1,1]$. The total number of iterations and errors for different polynomial degree approximations are shown in \cref{tabtestconv} and \cref{figconv1d}. We observe a consistent error decay as the approximation degree is increased, while the number of iterations of the policy iteration algorithm remains independent. The use of a monomial basis for approximating the value function yields similar results as using a set of orthogonal Legendre polynomials.

\begin{table}[!ht]
	\centering
	\begin{tabular}{ccccccc}
		\hline\\
		& \multicolumn{3}{c}{Monomial basis} & \multicolumn{3}{c}{ Legendre basis}\\
		\cmidrule(lr){2-4}\cmidrule(lr){5-7}\\
		$n$ &error $V_n$&error $u_n$& iterations &error $V_n$&error $u_n$& iterations\\
		\cmidrule(lr){1-1}\cmidrule(lr){2-2}\cmidrule(lr){3-3}\cmidrule(lr){4-4}\cmidrule(lr){5-5}\cmidrule(lr){6-6}\cmidrule(lr){7-7}\\
		2 & 1.4128 & 0.5919 &71 & 1.4035 & 0.5846 & 71\\
		4 & 0.3643 & 0.1468 &79 & 0.8184 & 0.3172 & 71\\
		6 & 0.0202 & 9.80$\times 10^{-4}$ &75& 0.0224 & 0.0101 & 75 \\
		8 & 6.39$\times 10^{-4}$ & 3.39$\times 10^{-4}$ &75&6.97$\times 10^{-4}$ &3.55$\times 10^{-4}$&75\\
		10 & 1.43$\times 10^{-5}$ &8.01$\times 10^{-6}$ &76 & 1.45$\times 10^{-5}$& 8.23$\times 10^{-6}$&77\\
		\hline
		\\
	\end{tabular}
	\caption{1D polynomial approximation for the HJI equation with nonquadratic running cost. The number $n$ denotes the polynomial degree of the approximated value function with a monomial and Legendre basis elements. Errors are shown for the value function $V_n(x)$ and for the optimal feedback $u_n(x)$. The number of iterations includes the total number of subiterations related to the disturbance update.}\label{tabtestconv}
\end{table}

\subsection{Robust control of nonlinear parabolic PDEs} In the following tests, we study the design of robust feedback controllers for a class of nonlinear parabolic PDEs including nonlinear advection and reaction terms such as the Burgers or Allen-Cahn equations. In order to embed the infinite-dimensional PDE dynamics into our finite-dimensional robust control framework, we semi-discretize the dynamics in space by means of a pseudospectral collocation method, leading to nonlinear dynamical systems with as many dimensions as degrees of freedom in the semi-discretization.
\subsubsection*{Test 1: a viscous Burgers' equation with forcing term} We study the robust optimal feedback stabilization of a control system related to a  Burgers' equation with a nonlinear reaction term. The control problem is stated as
\begin{align*}
\underset{u(\cdot)\in \cU}{\min}\underset{w(\cdot)\in \cW}{\max}\;&\cJ(u,w;X_0):=\int\limits_0^\infty \|X(\cdot,t)\|^2_{L^2(\cI)}+\|u(t)\|_R^2-\gamma^2\|w(t)\|_P^2\, dt\\
\text{subject to}&\\
(\Sigma_{\infty})&\left\{
\begin{array}{ll}
\partial_{t}X(\xi,t) &=\sigma \partial_{\xi\xi}X(\xi,t) +X(\xi,t)\partial_{\xi}X(\xi,t)
+ 1.5\, X(\xi,t)e^{-0.1X(\xi,t)}\\[1.5ex]&\qquad\ldots +\chi_{\omega_1}(\xi)w(t)
+\chi_{\omega_2}(\xi)u(t)\,,\\[1.5ex]
X(\xi_l,t)&=X(\xi_r,t)=0\,,\quad t\in \R^+,\\[1.5ex]
X(\xi,0)&=X_0(\xi)\,,\quad \xi\in \cI\,,
\end{array}
\right.
\end{align*}
where $(\xi,t) \in \mathcal{I}\times \mathbb{R}^+$ with
 $\cI=(-1,1)$ and $\cU=\cW=L^2([0,+\infty))$. Here $\chi_{\omega_1}$ and $\chi_{\omega_2}$ are indicator functions for the action of the disturbance and the control supported over subsets $\omega_1$ and $\omega_2$, respectively. The additional source term $1.5 X(\xi,t)e^{-0.1X(\xi,t)}$ prevents the origin to be origin asymptotically stable in the uncontrolled/unperturbed case $u=w=0$. In our numerical tests the viscosity parameter is set $\sigma=0.2$.

 In order to reduce the infinite dimensional state equation $\Sigma_{\infty}$ to a finite dimensional dynamical system, we apply a pseudospectral collocation method in space based on Chebyshev polynomials (see e.g. \cite{os06,qvbook}, Chapter 4) which leads to a state
 space representation of the  form
 \begin{align}\label{eq:kk1}
 \dot X(t)=AX(t)+(X(t))\circ(DX(t))+1.5 X(t)\circ e^{-0.1X(t)}+Cw(t)+Bu(t),
 \end{align}
 where the discrete state $X(t)=(X_0(t),\ldots , X_d(t))^t\in \mathbb{R}^{d+1}$ corresponds to the approximation of $X(\xi,t)$ at $d$ collocation points, $X_i(t)=X(\xi_i,t)$ with
 $\xi_i=-cos(\pi i/(d+1))$, $i=1,\ldots,d$, while $X(t)\circ DX(t)$ and $X\circ e^{-0.1X}$ denote  coordinatewise operations. The matrices $A\in \mathcal{M}^{d\times d}$, $D\in \mathcal{M}^{d\times d}$, $B\in \mathbb{R}^d$ and $C\in \mathbb{R}^d$ are
 finite dimensional approximations of the Laplacian, gradient, control and disturbance operators, respectively.  The $L^2$ norm of the state penalization of the running cost is approximated accordingly. Note that such a state-space representation is consistent with \cref{as:1}, as for each $i=1,\ldots, d$ the state equation reduces to
 \begin{align*}
 \dot X_i(t)&=A_{i,1}X_1(t)+\hdots+A_{i,d}X_d(t)+X_1(t)D_{1,d}X_1(t)+\hdots+X_d(t)D_{i,d}X_d(t)\\
 &\qquad+(1.5X_i)e^{-0.1X_i}+C_iw(t)+B_iu(t),
 \end{align*}
 where $A_{i,\cdot}$ denotes the $i$-th row of the matrix $A$. In this case, the dynamics have a separability degree $n_f= 2d+1$.  We shall solve the HJI equation
 \begin{align}
 DV_\gamma(x)^tf(x)+ \frac14 DV_\gamma(x)^t\left(\frac{1}{2\gamma^2}h(x)P^{-1}h(x)^t-g(x)R^{-1}g(x)^t\right)DV_\gamma(x)+\ell(x)=0 \,,
 \end{align}
 with $f(x),h(x),g(x)$ and $\ell(x)$ generated with the aforedescribed semi-discretization \eqref{eq:kk1} . The value function $V(x)$ is approximated with $d=12$ over $\Omega=(-2,2)^{12}$  with a monomial basis of total degree up to 4.

In the following numerical tests we will study the effect of the location of domains
$\omega_1$ and $\omega_2$ on the solution of the HJI problem and in particular on the resulting $\gamma^*$. We further investigate numerically the behaviour of the HJI control in the case $\gamma$ is chosen differently from the optimal $\gamma^*$, and we compare the HJB and HJI controllers in terms of their robustness against additive noise in the state equation.

 \noindent
 {\bf{Effect of $\omega_1$ and $\omega_2$.}}
 Here we are interested in the effect of the control and disturbance supports on the solution of the HJI equation. In particular, we concentrate  on the consequences in terms of the $\cH_{\infty}$-norm $\gamma^*$, the smallest value for which the HJI equation leads to an asymptotically stabilizing feedback control law. This value is computed by a bisection algorithm discarding values of $\gamma$ where \cref{alg:sga1} fails to converge.
For this test, we  consider the dynamics \eqref{eq:kk1} without the exponential source term. \cref{table:5} reports three cases varying $\omega_1$ and $\omega_2$.   For the largest
ratio of size of control to size of disturbance, we obtain the smallest $\gamma^*$, see case 1.
Case 3 arises from case 1 by increasing the support of the disturbance while decreasing the support of the control,  leading to a significant increase of $\gamma^*$. In case 2 we decrease the disturbance support compared to case 3, in such a way that $w_1\cap w_2=\emptyset$, which results in decrease of $\gamma^*$. The three cases are consistent with the expectation that as the ratio $|\omega_2|/|\omega_1|$ increases, stabilization and noise cancellation capabilities are improved resulting in a smaller $\gamma^*$ value.

\begin{figure}[!ht]
	\centering
	\subfloat[]{\includegraphics[width=0.5\textwidth]{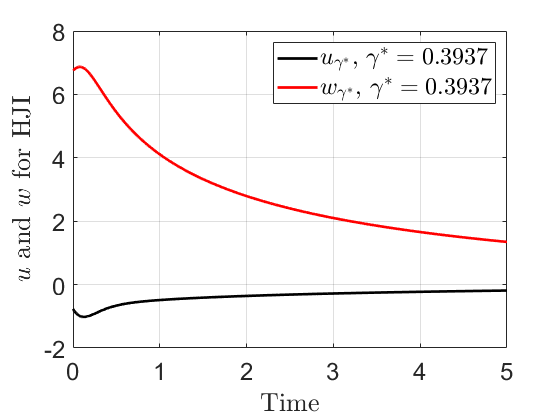}}\hfill
	\subfloat[]{\includegraphics[width=0.5\textwidth]{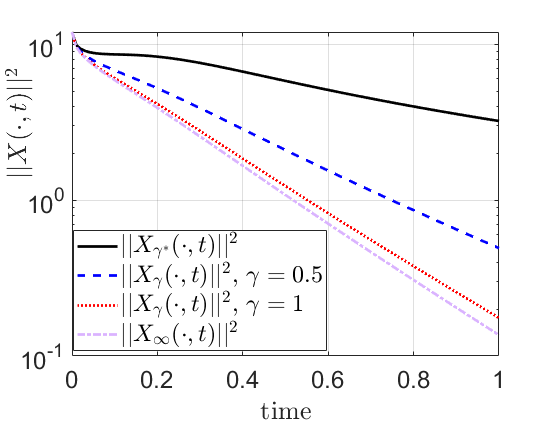}}\hfill
	\caption{Test 1 (viscous Burgers' equation without exponential source term), influence of $\gamma$.  Parameters $R=0.1$ and $P=1$, initial condition: $X(\xi,0)= sign(\xi)$. (a) Disturbance $w_{\gamma^*}$ and control $u_{\gamma^*}$ for HJI. (b) Behaviour of $\norm{X(\xi,t)}^2_{L^2(\cI)}$  as $\gamma$ increases.}
	\label{fig:figure1}
\end{figure}

\begin{table}[!ht]
\centering
 \begin{tabular}{c c c c c}
  \hline
   & support of disturbance ($\omega_1$) & support of control ($\omega_2$) & $\gamma^*$
 \\  \hline \hline
  case 1 & $(0.5,0.8)$ &    $(-1,1)$      & 0.3937  \\ \hline
  case 2 & $(-1,-0.5)$ & $(0.5,0.8)$ & 0.8087 \\ \hline
  case 3 & $(-1,1)$ & $(0.5,0.8)$ & 3.1281\\
  \hline
 \end{tabular}
 \caption{Disturbance support $\omega_1$, control support $\omega_2$, and  $\gamma^*$ for $R=0.1$ and $P=1$.}  \label{table:5}
\end{table}
\noindent
 {\bf{Effect of choosing $\gamma$}.}  For the following test, $\omega_1$ and $\omega_2$ are chosen as in case 1, above. We report that similar qualitative results were obtained in the other cases. We begin by setting $\gamma=\gamma^*$, results for  the initial condition $X(\xi,0)= sign(\xi)$ are shown in \cref{fig:figure1} (a) for $u$ and $w$ given by the feedback formulas
\begin{align}
u_{\gamma}(x)=-\frac12 R^{-1}g(x)^tD V_\gamma(x)\,,\qquad w_{\gamma}(x)=\frac{1}{2\gamma^2} P^{-1}h(x)^tD V_\gamma(x)\,.
\end{align}
For the results in  \cref{fig:figure1}(b) the HJI equation is solved for an  increasing sequence of values for $\gamma$. We observe that the corresponding norms of the trajectories  $\norm{X(\xi,t)}^2_{L^2(\cI)}$ converge as $\gamma$ increases to the graph for $\norm{X(\xi,t)}^2_{L^2(\cI)}$ which corresponds to the  HJB synthesis, equivalent to $\gamma=+\infty$. We report that	 for $\gamma^*$ the graph of  $ \norm{X(\xi,t)}^2_{L^2(\cI)}$ also tends to  0 as $t\to \infty$, but at a much smaller rate.
\noindent
\begin{figure}[!ht]
\centering
\subfloat[]{\includegraphics[width=0.5\textwidth]{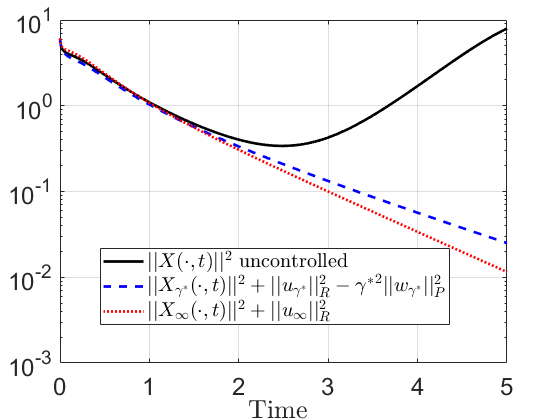}}\hfill
\subfloat[]{\includegraphics[width=0.5\textwidth]{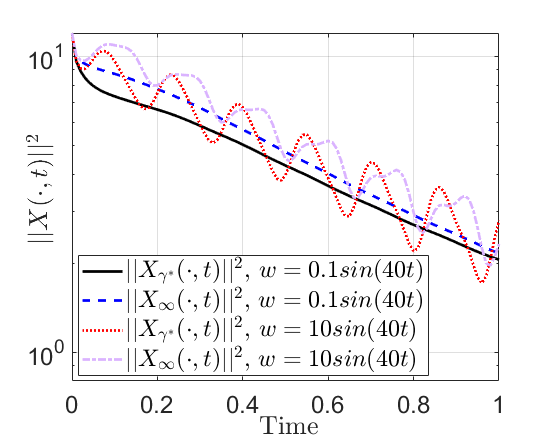}}\hfill
\caption{Test 1 (viscous Burgers' equation with exponential source term), $R=0.01,P=1$, $\gamma^*=0.125$. $X(\xi,0)= sign(\xi)$,  and sinusoidal disturbance. 
(a) Evolution of the running cost $\norm{X(\xi,t)}^2_{L^2(\cI)}+\norm{u(t)}^2_{R}-{\gamma^*}^2\norm{w(t)}^2_{P}$ with optimal feedback and disturbance. (b) Evolution of the state penalization $\norm{X(\xi,t)}^2_{L^2(\cI)}$ for HJB and HJI feedback laws under different sinusoidal disturbances.}
\label{fig:noise1}
\end{figure}

\noindent
{\bf{Comparison between HJI and HJB feedback laws in the presence of noise.}}
We compare the behaviour of the \eqref{eq:kk1} in the case that the control is computed by means of the HJI eqn. with optimal $\gamma^*$, and the HJB feedback law corresponding $\gamma=+\infty$. Different types of noise signals $w$ are tested. Here we include the exponential source term in the Burgers' equation, we set $R=0.01$ as the weight in the running cost for the control, and $\omega_1=\omega_2=(-0.8,0.5)$.

With the additional source term, the uncontrolled-unperturbed system is unstable. This is highlighted in \cref{fig:noise1} {(a)}, where the norm of the state of the uncontrolled system with initial condition $X(\xi,0)= sign(\xi)$ is compared against the state norms produced by the  HJB and HJI feedback laws with $w=0$. We note that both feedback laws generate an asymptotically stable closed loop system.

Next we  compare  the HJI/HJB synthesis capabilities to reject additive  noise. We compute $V_{\gamma^*}$ and $V_{\infty}$ to generate the  associated feedback control laws, which are evaluated for different choices of the perturbation signal $w(t)$. For a sinusoidal noise of the form  $w= \eta \sin(\omega t)$,  the total state contributions $\norm{X}^2 = \int_0^\infty \|X(t)\|^2\,dt$ as well as the total costs for HJI and HJB closed loops, $\cJ(u_{\gamma},w=0;X_0)$ and $\cJ(u_{\infty},w=0;X_0)$ respectively, are presented in \cref{table:noise1} for $X(\xi,0)= sign(\xi)$.
\noindent
\begin{table}[!ht]
	\centering
	\begin{tabular}{c c c c c c}
		\hline
	& $w(t)$& $\norm{X}^2$ with $u_{\gamma}$ & $\norm{X}^2$ with $u_{\infty}$	& $\cJ(u_{\gamma},0;X_0)$ &$\cJ(u_{\infty},0;X_0)$
		\\  \hline \hline	
	& $0.1sin(10t)$   & 2.371  & 2.666       & 2.377     & 2.669  \\ \hline
	& $10sin(10t)$    & 3.891  & 4.278       &4.184       & 4.605 \\ \hline
	 & $0.1sin(40t)$ & 2.361   & 2.659        & 2.369    & 2.661 \\ \hline	
	  & $10sin(40t)$  & 2.676  & 3.022       & 2.837       & 3.149    \\ \hline
		
	\end{tabular}
	\caption{HJI and HJB feedback laws under sinusoidal disturbances, Burgers' equation with exponential source term, $\cJ(u,w=0;X_0)=\int_0^\infty( \|X(t)\|^2+\norm{u}^2_{R})\,dt$.}\label{table:noise1}
\end{table}

We note that the state cost given by the HJI synthesis is consistently smaller than associated to the HJB feedback control law.  This can also  be confirmed in  \cref{fig:noise1} { (b)}. We proceed with other choices of the perturbation signals given by:
\begin{itemize}
	\item $w_1(t)=\kappa w_{\gamma}+\eta sin(\omega t)$
	\item $w_2(t)=\kappa u_{\gamma}+\eta sin(\omega t)$
	\item $w_3(t)=\kappa u_{\infty}+\eta sin(\omega t)$,
\end{itemize}
where $\kappa$ and $\eta$ are parameters, and  $w_{\gamma}$, $u_{\gamma}$ and $u_{\infty}$ are the HJI disturbance, the HJI control,  and the  HJB control signals respectively, the latter two obtained by simulating first the unperturbed dynamics.
\begin{figure}[!ht]
	\centering
	\subfloat[]{\includegraphics[width=0.5\textwidth]{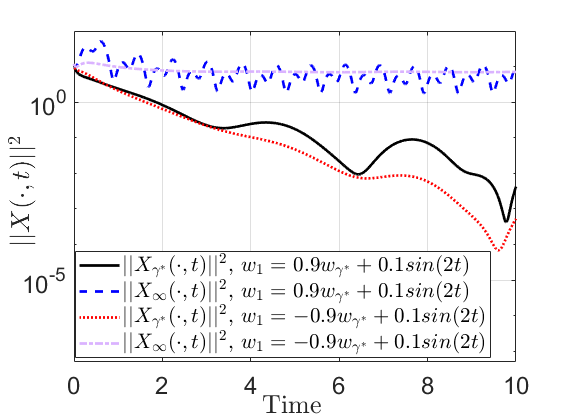}}\hfill	
	\subfloat[]{\includegraphics[width=0.5\textwidth]{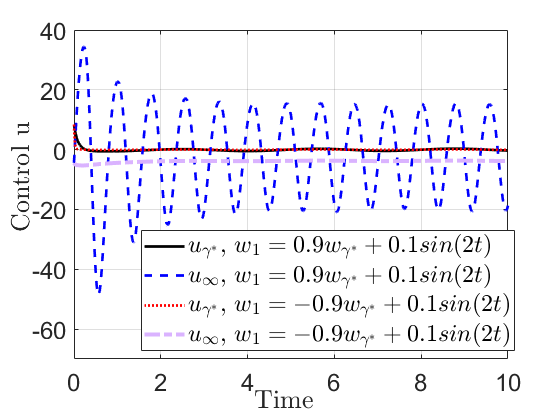}}\\
	\subfloat[]{\includegraphics[width=0.5\textwidth]{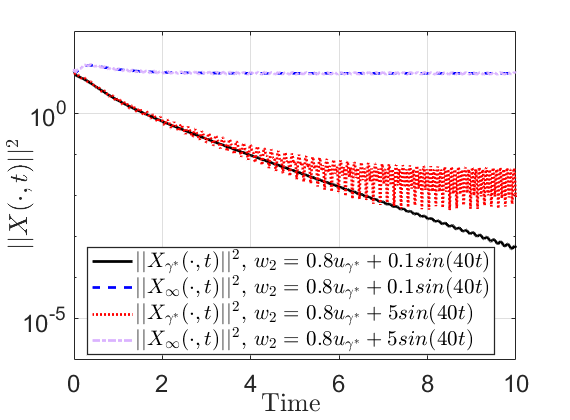}}\hfill
	\subfloat[]{\includegraphics[width=0.5\textwidth]{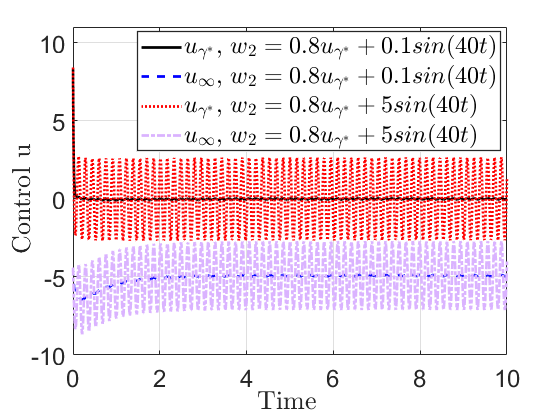}}\\
	\subfloat[]{\includegraphics[width=0.5\textwidth]{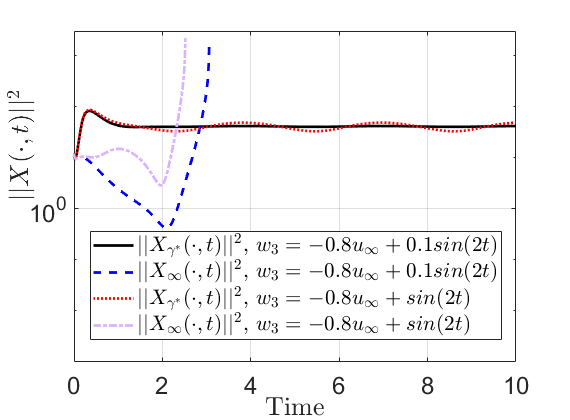}}\hfill
	\subfloat[]{\includegraphics[width=0.5\textwidth]{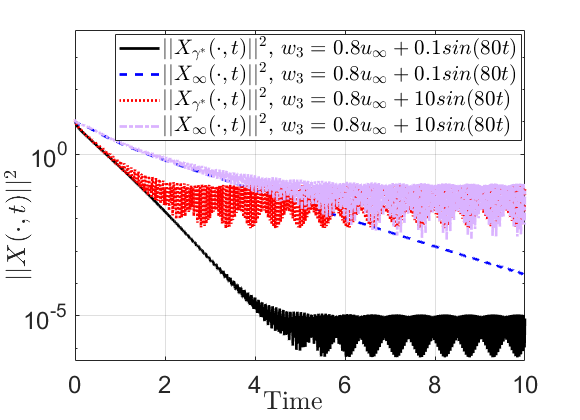}}
	\caption{Test 1(viscous Burgers' equation with exponential source term) comparing HJI and HJB feedback laws under different disturbances,  $R=0.01$, $P=1$, $X(\xi,0)=sign(\xi)$, and $\gamma^*=0.125$. }
\label{fig:noise2}
\end{figure}

The results {(a)-(b)} in the first row of \cref{fig:noise2}, are obtained for $w_1$, for the values of $\kappa$ and $\eta$  indicated in the figures.
Unlike the results with sinusoidal perturbation only shown in \cref{table:noise1}, where there is no significant difference between the HJI and the HJB closed loops, here the HJI synthesis asymptotically stabilizes the dynamics to the origin, whereas the HJB feedback law fails.
Also, from \cref{fig:noise2} {(c)} and (d) we observe that in the presence of the perturbation $w_2(t)$, the HJI feedback controls achieve stabilization to the origin, while the HJB feedback fails in this respect.
In \cref{fig:noise2} {(e)} and (f) we present results for the disturbance $w_3(t)$ with different values of $\kappa$. For $\kappa=-0.8$ the HJB synthesis not only fails to stabilize but generates a finite time blow up, while the HJI feedback is again successful. For $\kappa = 0.8$ both control strategies reject the disturbance, with the HJI law having a more effective transient behaviour.

\subsubsection*{Test 2: the Degenerate Zeldovich equation} In this second test we consider the stabilization of the following diffusion-reaction model with Neumann boundary conditions arising in combustion theory \cite{GK04}
\begin{align*}
\partial_{t}X(\xi,t) &=\sigma \partial_{\xi\xi}X(\xi,t)+X(\xi,t)^2-X(\xi,t)^3+\chi_{\omega_1}(\xi)w(t)+\chi_{\omega_2}(\xi)u(t)\,,\\
\partial_{\xi}X(\xi_l,t)&=\partial_{\xi}X(\xi_r,t)=0\,, t\in \R^+\,,\\
X(\xi,0)&=X_0(\xi)\,,\quad \xi\in \cI\,,
\end{align*}
with $(\xi,t)\text{ in}\;\cI \times\R^+$, $\cI=(-1,1)$, $\sigma=0.5$, $\omega_1=\omega_2=(-0.8,0.5)$. In this case $X\equiv 0$ is an unstable and $X\equiv 1$ is a stable steady state for the uncontrolled state equation. The dynamics are semi-dscretized in space with the same collocation method as in Test 1. The reduced state space is chosen to be $\Omega=(-2,2)^{12}$, and the basis functions for the value function  are monomials up to degree 4. We take $R=0.01$, $P=1$. Then the $\cH_{\infty}$-norm is found to be $\gamma^*=0.125$.
\begin{figure}[!ht]
	\centering
	\subfloat[]{\includegraphics[width=0.5\textwidth]{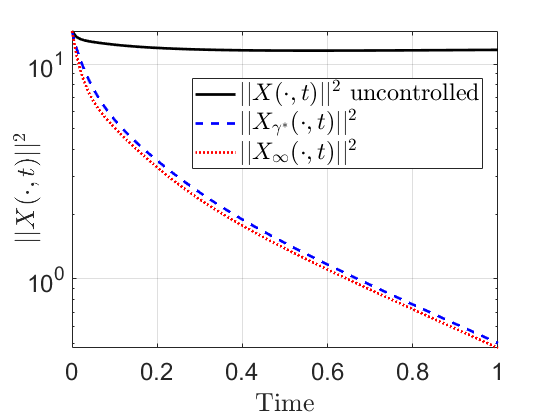}}\hfill
    \subfloat[]{\includegraphics[width=0.5\textwidth]{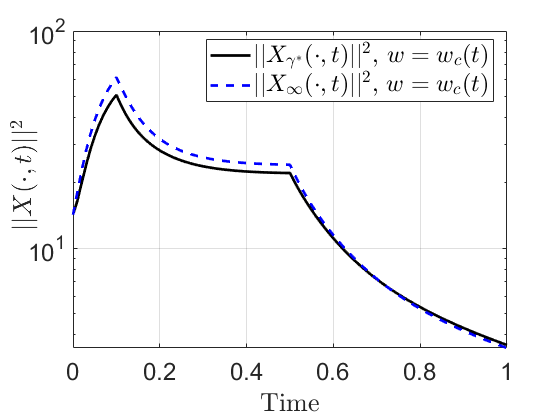}}\hfill
	\caption{Test 2: Degenerate Zeldovich equation, $R=0.01$, $P=1$, and $X(\xi,0)=cos(2\pi\xi)sin(\pi\xi)+1$.
		{(a)} Evolution of $\norm{X(\xi,t)}^2_{L^2(\cI)}$ for uncontrolled, HJI and HJB feedback laws. {(b)}
Evolution of $\norm{X(\xi,t)}^2_{L^2(\cI)}$ with the piecewise constant disturbance \eqref{eq:wc}. Total
		HJB cost=10.98, total HJI cost=10.20.}
	\label{fig:zelnoise1}
\end{figure}
In \cref{fig:zelnoise1} {(a)}, we depict the norm of the uncontrolled solution and compare it to the HJI and HJB controlled solutions. We observe the two feedback laws generating a similar response in the absence of noise. We also compared HJI and HJB with sinusoidal perturbations and found only marginal differences between the two cases, with a tendency or HJI to be superior in the case of small frequencies and vice versa for higher frequencies.

As in Test 1, we continue to compare HJI to HJB for non-sinusoidal disturbances, and we first consider the piecewise constant signal given by
\begin{align}
w_c(t)=
\begin{cases}
30 \quad & 0\leq t<0.1\\
10 \quad & 0.1\leq t\leq 0.5\\
0.5 \quad & t>0.5.
\end{cases}\label{eq:wc}
\end{align}
For this case we see from \cref{fig:zelnoise1} {(ii)} that in the transient phase the HJI synthesis is superior to the HJB closed-loop.

\begin{figure}[!ht]
	\centering
	\subfloat[]{\includegraphics[width=0.5\textwidth]{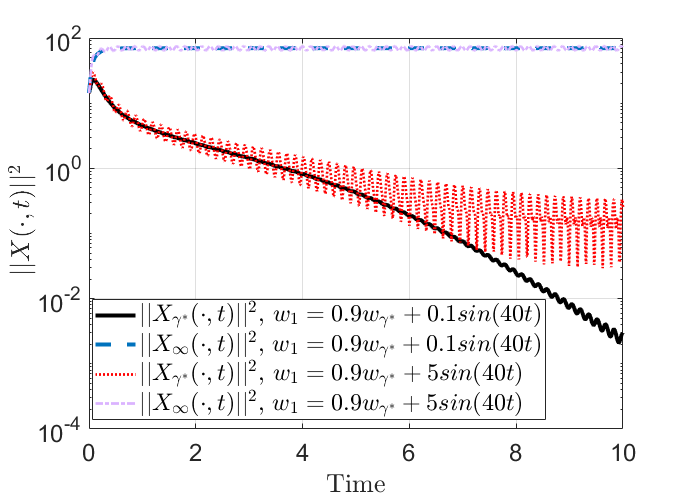}}\hfill
	\subfloat[]{\includegraphics[width=0.5\textwidth]{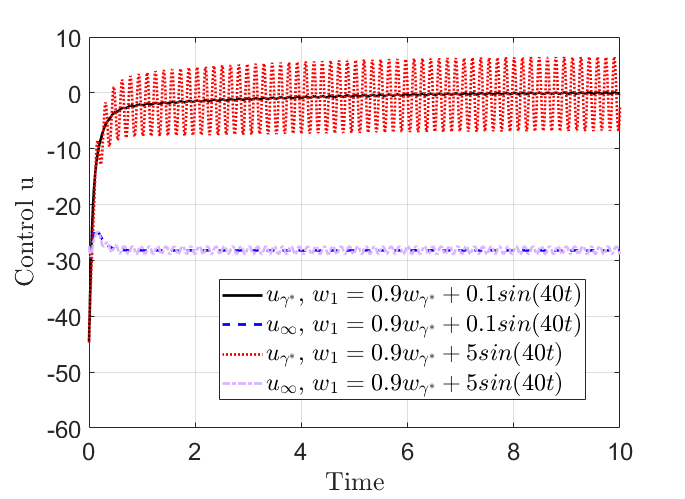}}\hfill
	\caption{Test 2: Degenerate Zeldovich equation, $R=0.01$, $P=1$, and $X(\xi,0)=cos(2\pi\xi)sin(\pi\xi)+1$. Disturbance: $ w(t)=0.9\, w_{HJI}+\eta sin(\omega t)$, {(a)} Evolution of $\norm{X(\xi,t)}^2_{L^2(\cI)}$, {(b)} Control signals for $u(t)$ in (a).}
	\label{fig:zelnoise2}
\end{figure}
Now we modify the disturbance term in the dynamics to observe the corresponding effect close to the  worst case scenario and choose
$w=0.9\, w_{\gamma^*}+\eta sin(\omega t)$.
We can see from \cref{fig:zelnoise2} that the HJB feedback control is not able to attenuate this disturbance level. The HJB feedback controller fails to steer the trajectories towards the unstable equilibrium, and the closed-loop dynamics remain attracted to $X=1$. In the presence of the same disturbance, the HJI controller succeeds in directing the state to the desired equilibrium steers the state to the origin. The choice of signals $w_2$ and $w_3$ as in Test 1 leads to a similar type of behaviour.

\subsubsection*{Test 3: Stabilization under model  uncertainty}
We conclude by studying the case of state-dependent perturbations. Alternatively, we can interpret this setting as modelling  parametric uncertainties in the dynamics as for instance, uncertainties in the diffusion and/or reaction coefficients. We first investigate the dynamics
\begin{align*}
\partial_{t}X(\xi,t) &=\sigma \partial_{\xi\xi}X(\xi,t)-X(\xi,t)^3+X(\xi,t)w(t)+\chi_{\omega_2}(\xi)u(t)\,,\\
\partial_{\xi}X(\xi_l,t)&=\partial_{\xi}X(\xi_r,t)=0\,,\quad t\in \R^+\,,\\
X(\xi,0)&=\kappa(\xi-1)^2(\xi+1)^2\,,\quad \xi\in \cI\,,
\end{align*}
with $(\xi,t)\text{ in}\;\cI \times\R^+$, $\cI=(-1,1)$, $\omega_2=(-0.8,0.5)$ and $\kappa$ to be specified below. As in the previous tests, the computational domain is $\Omega=(-2,2)^{12}$.

We compare between the HJB and HJI synthesis for trajectories generated by choosing the constant perturbation $w=1$ in the above model. Note that for $w=1$ the above reaction-diffusion model corresponds to the  Newell-Whitehead /Allen-Cahn equations
\begin{equation*}
\partial_{t}X(\xi,t)=\sigma \partial_{\xi\xi}X(\xi,t)-X(\xi,t)^3+X(\xi,t)+\chi_{\omega_2}(\xi)u(t)
\end{equation*}
 The HJB synthesis does not incorporate the destabilizing term $X(\xi,t)$ . We study the behaviour of the closed loop under additional perturbations of the type $\eta sin(\omega t)$.

\begin{figure}[!ht]
	\centering
	\subfloat[]{\includegraphics[width=0.5\textwidth]{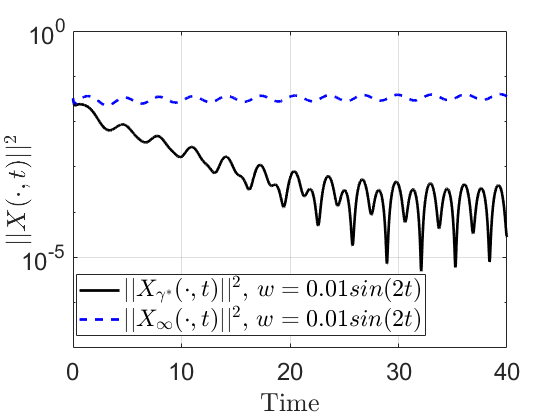}}\hfill
	\subfloat[]{\includegraphics[width=0.5\textwidth]{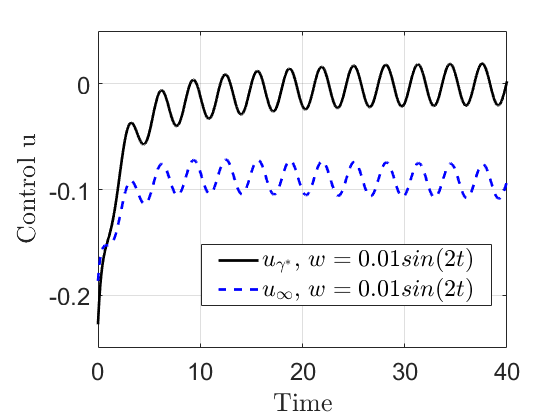}}\hfill
	\subfloat[]{\includegraphics[width=0.5\textwidth]{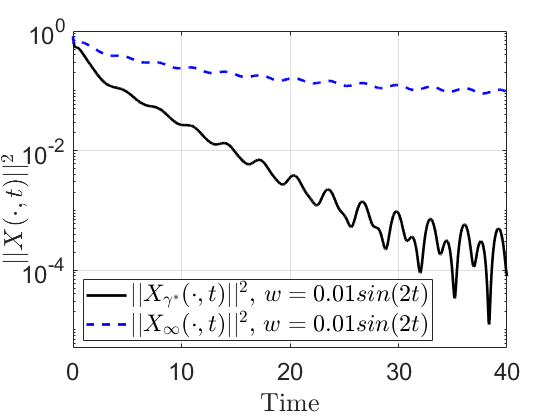}}\hfill
	\subfloat[]{\includegraphics[width=0.5\textwidth]{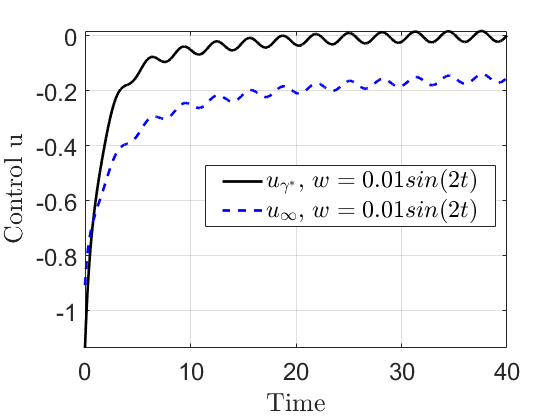}}\hfill
	\subfloat[]{\includegraphics[width=0.5\textwidth]{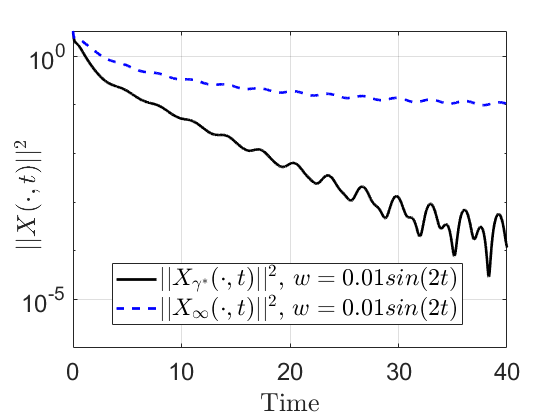}}\hfill
	\subfloat[]{\includegraphics[width=0.5\textwidth]{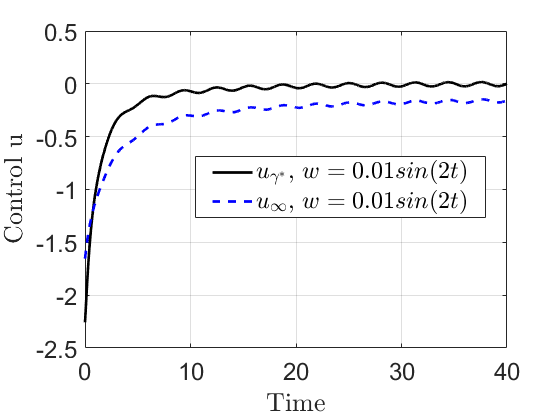}}\hfill
	\caption{Test 3, state-dependent disturbance in the Newell-Whitehead equation, $R=0.5$, $P=1$, $\sigma=0.5$, $\omega_2=(-0.3,0.5)$, $\gamma^*=1.63$ and $X(\xi,0)=\kappa(\xi-1)^2(\xi+1)^2$. State cost evolution and control signals for: (a)-(b) $\kappa=.1$, (c)-(d) $\kappa=.5$, (d)-(e) $\kappa=1$.}
	\label{fig:statenoise}
\end{figure}

In \cref{fig:statenoise} we close the loop with HJI and HJB based control, where in the three consecutive rows we increase the magnitude of the initial condition by considering the cases $\kappa=0.1$,  $\kappa=0.5$, and $\kappa=1$, respectively.

For all the three cases the HJI controller stabilizes the system towards the origin and the three corresponding controls tend to $0$ as time increases. The HJB feedback law, which does not anticipate the modelling error in the reaction term, fails to stabilize this system to $0$.

Finally, we consider a model uncertainty for a cubic reaction term in
\begin{align*}
\partial_{t}X(\xi,t) &=\sigma \partial_{\xi\xi}X(\xi,t)+X(\xi,t)^3w(t)+\chi_{\omega_2}(\xi)u(t)\,,\qquad\text{in}\;\cI\times\R^+\,,\\
\partial_{\xi}X(\xi_l,t)&=\partial_{\xi}X(\xi_r,t)=0\,,\quad t\in \R^+\,.
\end{align*}
\begin{figure}[!ht]
	\centering
    \subfloat[]{\includegraphics[width=0.33\textwidth]{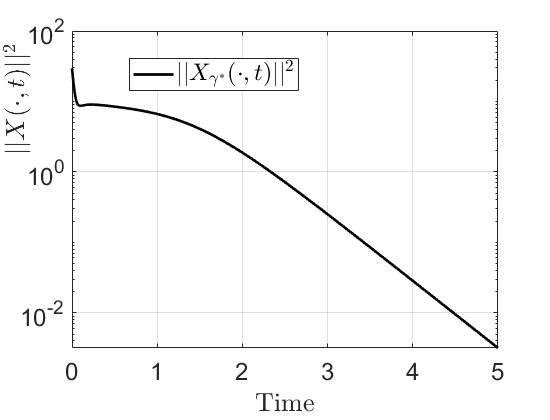}}\hfill
	\subfloat[]{\includegraphics[width=0.33\textwidth]{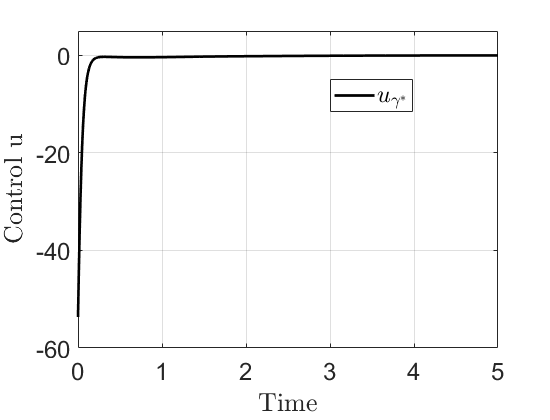}}\hfill
	\subfloat[]{\includegraphics[width=0.33\textwidth]{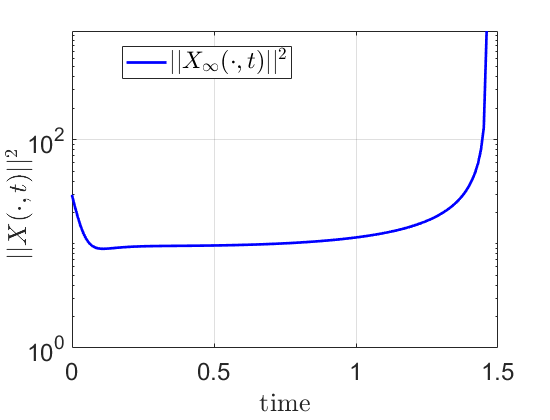}}\hfill
	\caption{Test 3, $X^3$ state-dependent disturbance in the Newell-Whitehead equation, $R=0.01$, $P=1$, $\sigma=0.5$, $\omega_2=(-0.8,0.5)$, $\gamma^*=8.22$ and $X(\xi,0)=3(\xi-1)^2(\xi+1)^2$. {(a)} Evolution of  $\norm{X(\xi,t)}^2_{L^2(\cI)}$ with HJI feedback law, {(b)} Optimal HJI control signal for (a) and,
				{(c)} Evolution of $\norm{X(\xi,t)}^2_{L^2(\cI)}$ with HJB feedback (the closed-loop becomes unstable).}
	\label{fig:statecubicnoise}
\end{figure}
In this case a robust control design is critical, as the cubic reaction term can generate a finite time blow-up of the dynamics. We observe that for the case $w=1$, the HJB control fails to stabilize the dynamics whilst the HJI closed-loop is asymptotically stable.  This is documented in \cref{fig:statecubicnoise}

\subsection*{Remarks on CPU  times} \cref{tab:cpu} provides CPU-assembly and CPU-iterative times for solving the  HJI equation related to Tests 1 and 2. Tests were run on a muti-core architecture 8x Intel Xeon E7-4870 with 2,4Ghz, 1 TB of RAM. MATLAB pseudoparallelization distributes the tasks among 12 workers.
\begin{table}[!ht]
	\centering
	\setlength{\tabcolsep}{1mm}
	\begin{tabular}{cccc}
		\hline\\
		Test & Source term &CPU-assembly&CPU-iterative \\
		\cmidrule(lr){1-1}\cmidrule(lr){2-2}\cmidrule(lr){3-3}\cmidrule(lr){4-4}\\
		1 & (Without source term) & $1.161\times 10^4[s]$ & $3.5078\times 10^4[s]$\\
		1 & (With source term) &$2.1508\times 10^4[s]$ & $1.1808\times 10^5[s]$\\
		2 &                    & $3.0262\times 10^3[s]$ & $4.8975\times 10^4[s]$\\
		\hline
		\\
	\end{tabular}
	\caption{ CPU-assembly corresponds to the amount of time spent in offline assembly of the different terms of the Galerkin residual equation \eqref{gHJIg}. CPU-iterative refers to the amount of time spent inside \cref{alg:sga1}.}\label{tab:cpu}
\end{table}

\section*{Concluding remarks}
We have proposed a numerical scheme for the design of robust optimal stabilizing feedback controllers, and assessed its capabilities in the context of stabilization of nonlinear parabolic PDEs.  The overall technique consists of two main steps: i) the projection of the infinite-dimensional nonlinear parabolic dynamics onto a finite-dimensional relevant subspace by means of a pseudospectral collocation method  and ii) the synthesis of a robust control  by solving a Hamilton-Jacobi-Isaacs equation, for which we have resorted to a Galerkin-type method using a global polynomial ansatz for the value function, together with iterative techniques for the solution of the nonlinear HJI.
Overall, the proposed methodology proves to be successful in achieving the fundamental features of the $\cH_\infty$ control design, that is, to generate an optimal stabilization of the closed-loop with enhanced disturbance mitigation properties. We have numerically assessed that our synthesis scheme offers a more robust control law than the classical HJB/$\cH_2$ feedback. Furthermore, the proposed scheme mitigates the curse of dimensionality,  allowing to approximate the  $\cH_{\infty}$ synthesis for nonlinear  dynamics  of high dimensions. The developed framework also opens the door to the study of nonlinear feedback control of PDEs under uncertainty.

\section*{Acknowledgments}  The research of this paper was supported by the ERC advanced grant 668998 (OCLOC) under the EU's H2020 research program and the Imperial College London Research Fellowship program.
\bibliographystyle{siamplain}

\end{document}